\documentclass
{ifacconf}

\counterwithin*{section}{part}

\usepackage{graphicx}      
\usepackage[comma]{natbib}        

\usepackage{amssymb,dsfont}
\usepackage{amsmath,amsfonts}
\usepackage{tikz}
\usepackage{pgfplots}
\usetikzlibrary{intersections, pgfplots.fillbetween}
\pgfplotsset{compat=1.18} 
\usepackage{todonotes}
\usepackage{bbold}

\usepackage{cleveref}
\Crefformat{figure}{#2Fig.~#1#3}
\Crefmultiformat{figure}{Figs.~#2#1#3}{ and~#2#1#3}{, #2#1#3}{ and~#2#1#3}
\usepackage{caption}
\usepackage{subcaption}
\usepackage{graphicx}

\newcommand{\setdef}[2]{\left\{ \ #1\ \left|\ \vphantom{#1} #2\ \right.\right\}}
\newcommand{\N}{\mathds{N}}
\newcommand{\R}{\mathds{R}}
\newcommand{\E}{\mathds{E}}

\newcommand{\Rp}{\R_{\geq0}}
\newcommand{\vp}{\varphi}
\newcommand{\diag}{\text{diag}}

\begin{document}
\begin{frontmatter}

\title{Reinforcement Learning for Docking Maneuvers with Prescribed Performance\thanksref{footnoteinfo}} 

\thanks[footnoteinfo]{
This research has been conducted within the project frame of SeRANIS Seamless Radio Access Networks in the Internet of Space. The project is funded by dtec.bw Digitalization and Technology Research Center of the Bundeswehr, grant number 150009910. dtec.bw is funded by the European Union-Next Generation EU.
L.~Lanza and K.~Worthmann gratefully acknowledge funding by the Deutsche Forschungsgemeinschaft (DFG, German Research Foundation) -- {Project-IDs 471539468 and 507037103}. L.~Lanza further is grateful for funding by the Carl Zeiss Foundation (VerneDCt -- Project-ID 2011640173).
K.~Lux-Gottschalk gratefully acknowledges funding by the Ir\`{e}ne Curie fellowship.
}

\author[First]{Simon Gottschalk} 
\author[Second]{Lukas Lanza} 
\author[Second]{Karl Worthmann}
\author[Third]{Kerstin Lux-Gottschalk}

\address[First]{Universit\"at der Bundeswehr M\"unchen, Department of Aerospace Engineering, Germany(e-mail: simon.gottschalk@unibw.de).}
\address[Second]{Technische Universit\"at Ilmenau, Institute of Mathematics, 
Germany (e-mail: \{lukas.lanza,karl.worthmann\}@tu-ilmenau.de)}
\address[Third]{Eindhoven University of Technology, Department of Mathematics and Computer Science, Netherlands (e-mail: k.m.lux@tue.nl)}

\begin{abstract}                
We propose a two-component data-driven controller to safely perform docking maneuvers for satellites.
Reinforcement Learning is used to deduce an optimal control policy based on measurement data.
To safeguard the learning phase, an additional feedback law is implemented in the control unit, which guarantees the evolution of the system within predefined performance bounds.
We define safe and safety-critical areas to train the feedback controller based on actual measurements. 
To avoid chattering, a dwell-time activation scheme is implemented.
We provide numerical evidence for the performance of the proposed controller for a satellite docking maneuver with collision avoidance.
\end{abstract}


\end{frontmatter}
\section{Introduction} \label{sec:intro}
In recent years, traffic in space has significantly increased with new players sending satellites to space and old spacecrafts tumbling around as space debris. This adds further complexity to safely perform
collision-free docking maneuvers.
Furthermore, the ground station can only intervene time delayed. The complexity of the problem pushes classical optimal control approaches to their boundaries and current state-of-the art Machine Learning (ML) methods often lack performance guarantees. This motivates us to develop an automatic controller, which ensures successful and collision-free completion of docking maneuvers.

In general, path planning and collision avoidance is done by solving optimal control problems. Therefore, one defines an objective function and constraints, which include the equations of motion of the dynamical system as well as collision constraints. 
Typically, dynamic programming approaches are well suited to solve collision-constrained optimal control problems (e.g. \cite{Richter.2023}). Due to the high dimensionality of docking maneuvers in space, alternative approaches are needed. In the literature, other control strategies can be found, like the direct approach by \cite{Michael.2013} or the model predictive control idea, which is considered in \cite{Ravikumar.2020}.
However, the combination of high complexity, a large number of degrees of freedom, and variability w.r.t.\ initial conditions pushes these \textit{classical} optimal-control approaches to their limits. 
Recently, following the path of artificial intelligence, Reinforcement Learning (RL) approaches have emerged as a remedy, cf.~\cite{Bertsekas.2019}.
These data-based control approaches can handle high degrees of freedom and various scenarios without recomputing optimal control solutions.
For slight variations, it is sufficient to just execute the trained policy, e.g., a neural network, acting as a fast online controller.
However, so far, not much is known about performance guarantees for such control strategies, which is highly relevant for safety critical control tasks such as collision avoidance. 
 Thus, as an additional safeguarding mechanism, 
we use the \textit{funnel controller}
proposed in~\cite{ilchmann2002tracking}. 
Funnel control is a high-gain adaptive feedback controller with the following two advantages:
First, it achieves tracking of a given reference trajectory within predefined error margins.
Second, the tracking is achieved for unknown nonlinear multi-input multi-output systems. 
We highlight that no system equations are required; 
rather the following structural assumptions are made: 
well-defined relative degree, bounded-input bounded-output stable internal dynamics, and a high-gain property.
The latter means that the system can react sufficiently 
fast if only the input is excited by a large enough signal. 
%
Considering funnel control (or related concepts) as a safety filter has been topic of recent research, e.g., in combination with feedforward control in~\cite{drucker2023experimental}, robustifying model predictive control in~\cite{berger2023robust}, 
as switched controller in~\cite{bikas2024prescribed},
in the context of sampled-data control of continuous-time nonlinear systems in~\cite{schmitz2023safe},
and combining soft and hard constraints in~\cite{mehdifar2022funnel}, to name but a few.

In \cite{Saxena.2023_arxiv}, the authors present a RL approach, in which the design of rewards for deep Q-learning builds upon funnel functions for learning a control policy that enforces some Signal Temporal Logic specifications. 
In \cite{Xia.2023}, the authors present a robust RL control strategy that includes motion constraints for a vertical take-off and landing of an unmanned aerial vehicle on a moving target. 
In \cite{Berthier.2022} and \cite{Berthier.2021}, a trajectory tracking problem is presented for a floating satellite with commanded torques. We emphasize that the floating satellite in this manuscript surpasses the satellites model in its complexity. 

The main contribution of this publication is to add a robustness component to the training phase of a RL based controller such that the obtained policy guarantees to stay within a prescribed safety region. Previous work has already addressed the construction of a tracking controller with prescribed performance for nonlinear systems including a safeguard component in the learning process, cf.~\cite{lanza2023control,schmitz2023safe}. 
Here, we build upon these results and enhance the RL framework allowing for other control strategies than Q-learning as well. Moreover, we level its applicability to more complicated control systems such as the control of satellites in space, represented by a sophisticated satellite model.

This manuscript is structured as follows. 
In \Cref{sec:satelliteModel}, we model the satellite, which is equipped with a robot arm for docking maneuvers. 
Furthermore, we investigate its mathematical properties and verify the assumptions required to apply funnel control. 
In \Cref{sec:controlUnits}, we introduce the overall control objective, the RL policy, and the funnel controller. Then, we present our novel methodology combining RL and funnel control.
To avoid continuous action of the funnel controller, we implement an activation function and use a dwell-time scheme to prevent chattering.
In \Cref{sec:numerics}, we illustrate the two-component learning-based controller via a numerical simulation of a docking maneuver.

\ \\
\textbf{Notation.}
$\Rp:=[0,\infty)$. 
For $k\le m \le n \in \N$ and $v \in \R^n$, we denote with $v_{k:m}$ the components $(v_k,\ldots,v_m) \in \R^{m-k+1}$.
$W^{k,\infty}(I,\R^n)$ is the Sobolev space of all $k$-times weakly-differentiable functions
$f:I\to\R^n$ with $f,\dots, f^{(k)}$ essentially bounded. $0_{n\times m}$ for $m, n \in \N$ represents a zero matrix with $n$ rows and $m$ columns. Furthermore, with $\mathbb{1}_{n\times n}$ we denote an $n\times n$ identity matrix.
\section{Modeling the satellite} \label{sec:satelliteModel}
\begin{figure}[h]
	\begin{center}
		\includegraphics[width=0.4\textwidth]{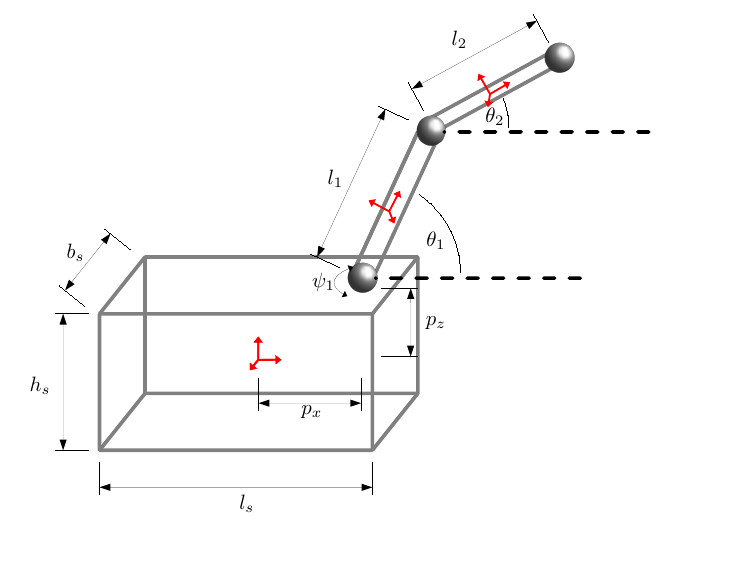}
  \vspace*{-2em}
		\caption{Schematic representation of the satellite model.}
		\label{fig:sketchSatellite}
	\end{center}
\end{figure}

In this section, we model the satellite as multibody system. We follow the steps from the book of \cite{Kortuem.1993}. The general structure of the satellite can be seen in~\Cref{fig:sketchSatellite}. It consists of a satellite body, which is described by the Cartesian coordinates of its center of mass $[x,y,z]$ as well as its roll $\phi$, pitch $\theta$ and yaw $\psi$ angle. Furthermore, attached to this body, we assume to have a robot arm. This arm consists of two elements. One is directly linked to the satellite's body. The spherical joint allows two degrees of freedom ($\psi_1$,$\theta_1$). The second element of the robot arm is attached to the first one by a revolute joint ($\theta_2$). We assume that we can manipulate the satellite by applying forces to accelerate the satellite body in $x$-, $y$- and $z$-direction, angular momentum to the satellite body and angular momentum to the robot arm.
In total, we have nine independent control inputs. The parameters of the satellite can be found in \Cref{tb:mearurements}.

\begin{table}[hb]
	\begin{center}
		\caption{Parameters of the Satellite Model}\label{tb:mearurements}
		\begin{tabular}{cccc}
			Variable & Unit & Value & Description \\\hline
			$m_1$ & [kg] & $300$ & mass of satellite body \\
			$m_2$ & [kg] & $1$ & mass of first robotic link \\
			$m_3$ & [kg] & $1$ & mass of second robotic link \\
			$l_s$ & [m] & $1$ & length of satellite body \\
			$b_s$ & [m] & $0.5$ & width of satellite body \\
			$h_s$ & [m] & $0.5$ & height of satellite body \\
			$r_1$ & [m] & $0.1$ & radius of first robotic link\\
			$r_2$ & [m] & $0.1$ & radius of second robotic link\\
			$l_1$ & [m] & $1$ & length of first robotic link\\
			$l_2$ & [m] & $1$ & length of second robotic link\\
			$p_x$ & [m] & $0.4$ & x-distance from center of gravity \\ &&&of satellite body to mount point\\
			$p_z$ & [m] & $0.25$ & z-distance from center of gravity \\ &&& of satellite body to mount point\\
		\end{tabular}
	\end{center}
\end{table}

For the derivation of the equations of motion, we introduce $q_1=[x,y,z,\phi,\theta,\psi]^\top$ as well as $q_2=[\theta_1,\psi_1,\theta_2]^\top$. Following the steps from \cite{Kortuem.1993}, the equations of motion read
\begin{align*}
	M(q_1(t),q_2(t))\ddot{q}_1(t)&=f(t,\dot{q}_1(t),\dot{q}_2(t),q_1(t),q_2(t)) \\
            & \quad +g(q_1(t),q_2(t))^\top u_{1:6}(t), \\
    \ddot{q}_2(t) &= u_{7:9}(t),
\end{align*}
where $u(t) \in \R^9$ represent the control inputs. Please note that, from now on, we will suppress the time dependencies of all quantities in this section for better readability.
We define the output
$y=(q_1,q_2)^\top$,
and may write the equations of motion in input/output form (invertibility given since $M$ is positive definite, see below) as
\begin{equation}  \label{eq:Model_input_output_form}
\begin{aligned} 
	\ddot{y} 
&=\begin{bmatrix}
 	M(y)^{-1}f(t,\dot y, y) \\ 0_{3x1}
 \end{bmatrix}  + \begin{bmatrix}
 	M(y)^{-1} g(y)^\top & 0_{3x3} \\  0_{3x3} & \mathbb{1}_{3x3}
 \end{bmatrix} u. 
\end{aligned}
\end{equation}
To check the assumption in order to apply the funnel theory, we need to check the positive definiteness of the matrix
\begin{align*}
   B(q_1,q_2):= \diag \left( M(q_1,q_2)^{-1} g(q_1,q_2)^\top, \mathbb{1}_{3x3}\right).
\end{align*}
Because of the block diagonal structure of this matrix, we can focus on the upper left block. Therefore, we focus on the mass matrix $M(q_1,q_2)$ and the matrix $g(q_1,q_2)$. Due to their long forms, we omit to write down the exact representation of the mass matrix $M$ and the matrix $g$. Instead, we introduce the center of mass $r_i(q_1,q_2)$ and the angular velocity $\omega_i(q_1,q_2)$ of body $i$. Together with the corresponding Jacobians
\begin{align*}
    J_i = \frac{\partial r_i(q_1,q_2)}{\partial q_1}, \quad 
    J_{i\omega} = \frac{\partial \omega_i(q_1,q_2)}{\partial q_1}, \quad \text{for } i=1,2,3,
\end{align*}
and the moments of inertia $I_{i}$ of body $i$, the mass matrix can be represented as (cf. \cite{Kortuem.1993}):
\begin{align*}
    M(q_1,q_2) = \sum_{i=1}^3 m_i J_i^\top J_i+ J_{i\omega}^\top I_{i} J_{i\omega}.
\end{align*}
Because of its structure, we can directly deduce that the mass matrix is symmetric and positive semi-definite and,
for $\theta\neq -\frac{\pi}{2}$, it is positive definite. We 
point out that all (reference) trajectories in the case study presented in \Cref{sec:numerics}
will not come close to the singularity $\theta= -\frac{\pi}{2}$. Since $M(q_1,q_2)$ is positive definite for all relevant states,
its inverse is positive definite as well. 
Similar, the matrix $g$ reads
\begin{align*}
    g(q_1,q_2)^\top &= \begin{bmatrix}
        J_1^\top A_z(\psi)A_y(\theta)A_x(\phi) \\
        J_{1\omega}^\top
    \end{bmatrix}
\end{align*}
with the rotation matrices:
\begin{align*}
    A_x(\phi) &= \small{ \begin{bmatrix}
        1 & 0 & 0 \\ 
        0 & \cos(\phi) & -\sin(\phi) \\
        0 & -\sin(\phi) & \cos(\phi) \\
    \end{bmatrix}}, 
    A_y(\theta) = \small{\begin{bmatrix}
        \cos(\theta) & 0 & -\sin(\theta) \\
        0 & 1 & 0 \\ 
        -\sin(\theta) & 0 & \cos(\theta) \\
    \end{bmatrix}}, \\
    A_z(\psi) &= \small{\begin{bmatrix}
        \cos(\psi) & -\sin(\psi) & 0\\
        -\sin(\psi) & \cos(\psi) & 0\\
        0 & 0 & 1 \\ 
    \end{bmatrix}}.
\end{align*}
The symbolic inversion of the mass matrix is computationally complex. Because of the bijection defined by $w = M^{-1}v$ for $v\in \R^6$ and the transformation
\begin{align*}
    & v^\top M^{-1}g^\top v = v^\top M^{-1}g^\top M^\top M^{-\top}v \\
    &= (M^{-\top}v)^\top g^\top M^\top (M^{-\top}v) = w^\top g^\top M^\top w 
    = w^\top g^\top M w,
\end{align*}
we focus on the positive definiteness of the matrix $g^\top M$. We calculated the eigenvalues of the matrix $g^\top M+(g^\top M)^\top$ evaluated in state space values of a pre-specified fine grid numerically for $\phi,\theta,\psi \in [-\frac{\pi}{8},\frac{\pi}{8}]$ and $\theta_1,\psi_1,\theta_2 \in [-\pi,\pi]$. The occurrence of only positive eigenvalues provides numerical evidence for its positive definiteness.
Hence, system~\eqref{eq:Model_input_output_form} has well-defined relative degree two in the set
\begin{align*}
    \mathbb{X} := \left\{\begin{bmatrix}
       q_1 \\ q_2
    \end{bmatrix} \in \R^9 | \phi,\theta,\psi \in \left[-\frac{\pi}{8},\frac{\pi}{8}\right], \theta_1,\psi_1,\theta_2 \in [-\pi,\pi]\right\}.
\end{align*}
From this, we directly deduce that the system has the high-gain property, cf.~\cite[Rem.~1.3]{berger2021funnel}, and moreover, it has trivial (and thus stable) internal dynamics in~$\mathbb{X}$.
Therefore, system~\eqref{eq:Model_input_output_form} is accessible for funnel control in the above indicated area.
This means that, starting in $\mathbb{X}$, funnel control can be used to force the system to evolve within this area, cf.~\eqref{eq:ErrorPerformance}.
\section{Control objective and controller} \label{sec:controlUnits}
In this section, we introduce the control objective and the two-component controller to achieve that objective.

\subsection{Control objective} \label{Sec:ControlObjective}
\begin{figure}[ht]
\begin{subfigure}[t]{0.4\linewidth}
\begin{tikzpicture}[>=stealth,scale=0.55]
\draw[->] (0,0) --node[below]{$t$} (6.3,0);
\draw[->] (0,0) -- node[rotate=90,left=1em,above=0.1]{$\|e\|$} (0,3) ;
\draw[scale=2, dashed, domain=0:3, samples=100 , variable=\x, black] plot ({\x}, {1.2*exp(-\x)+0.05});
\draw[scale=2, domain=0:3, samples=100 , variable=\x, black] plot ({\x}, {abs(cos(deg(0.5*\x))*cos(deg(3*\x))*exp(-0.7*\x)});
\draw (1.4,1.4) -- node[above=0.7em,right=0.1em]{\footnotesize{$1/\vp$}} (1.7,1.8);
\end{tikzpicture}
\subcaption{Evolution of the error~$e$ within boundary~$1/\vp$.}
\label{Fig:ErrorInFunnel}
    \end{subfigure}
    \qquad 
    \begin{subfigure}[t]{0.4\linewidth}
    \begin{tikzpicture}[>=stealth,scale=0.55]
\draw[, name path= axis,->] (0,0) --node[below]{$t$} (6.3,0);
\draw[->] (0,0) -- node[rotate=90,left=1em, above=0.3em]{$\|e_2\|/\vp$} (0,3) ;
\draw[name path=upper, scale=2, dashed, domain=0:3, samples=100 , variable=\x, black] plot ({\x}, {1.2*exp(-\x)+0.05});
\draw[name path=lower, scale=2, dashed, domain=0:3, samples=100 , variable=\x, black] plot ({\x}, {0.75*(1.2*exp(-\x)+0.05});
\tikzfillbetween[of=upper and lower]{red, opacity=0.1};
\tikzfillbetween[of=lower and axis]{green, opacity=0.1};
\draw[scale=2, domain=0:3, samples=100 , variable=\x, black] plot ({\x}, {abs(1.9*sin(deg(1*\x))*cos(deg(0.8*\x))*exp(-0.9*\x)});
\draw (1.4,1.4) -- node[above=0.7em,right=0.1em]{\footnotesize{$1/\vp$}} (1.7,1.8);
\end{tikzpicture}
\subcaption{Safe (green) and safety-critical (red) region for~$e_2$.}
\label{Fig:SafeAndSafetycritical}
\end{subfigure}
\caption{Tracking error and funnel boundary.}
\label{Fig:ErrorAndFunnel}
\end{figure}
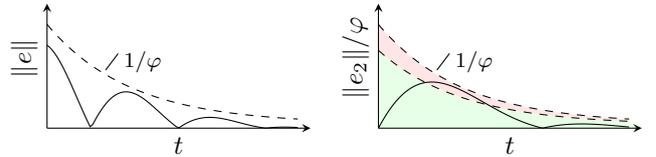

The aim is to use RL to derive an (optimal) control strategy from system data to perform a docking maneuver in space.
To safeguard the learning phase, 
the data-driven controller is equipped with 
an additional feedback controller (funnel control). 
This feedback controller is capable to compensate possible undesired control actions
such that the output $y = (q_1,q_2)$ of the system follows a given trajectory~$y_{\rm ref}$ with prescribed accuracy, i.e.,
\begin{equation} \label{eq:ErrorPerformance}
     \forall \, t \ge 0 \, : \ \| y(t) - y_{\rm ref}(t) \| < 1/\vp(t)
\end{equation}
for a user defined error tolerance~$1/\vp(t) > 0$.
This situation is illustrated in \Cref{Fig:ErrorInFunnel}.
Requirements on the reference~$y_{\rm ref}$ and the funnel boundary function~$\vp$ are presented in detail in \Cref{Sec:FeedbackLaw}.

\subsection{Reinforcement Learning}
As we know from Section \ref{sec:satelliteModel}, the satellite motion is described by $18$ variables ($q_1$, $q_2$, $\dot{q}_1$ and $\dot{q}_2$). This high complexity pushes classical optimal control approaches to their limits. 
However, AI approaches have shown that they can cope with high-dimensional problems. 
Thus, we apply the Proximal Policy Optimization (PPO) algorithm from~\cite{Schulman.2017}. In order to apply an RL approach, we need to define the underlying Markov Decision Process 
with the \textbf{State space~$S$} containing the states $[q_1^\top\, q_2^\top\, \dot{q}_1^\top\, \dot{q}_2^\top]^\top$ and the \textbf{Action space~$A$} representing the feasible control values, i.e., $[-0.75,0.75]^6\times [-0.15,0.15]^3$. The \textbf{transition probability} does not have to be specified explicitly. Trajectories are generated with the above introduced model. The initial position of the satellite is $s_0=0_{18 \times 1}$.
The \textbf{reward function~$r:S \times A \rightarrow \R$} will be defined in Section~\ref{sec:numerics}.

The idea of model-free online RL is based on the interaction between policy and environment in the form of a feedback loop. Typically, this control loop is defined for a discrete-time framework. Thus, we introduce the equidistant discretization points $t_0 < t_1 < \dots < t_N$ with $N \in \N$ and $\Delta h = t_{k+1}-t_k$ for all $k =1,\dots,N-1$. Thereon, the RL algorithm can be applied. 

PPO is a gradient based RL approach, which iteratively improves a parameterized policy $\pi_{\mu}:A \times S \rightarrow \R_+$ for parameters $\mu$. This randomized policy $\pi_{\mu}(\cdot,s)$ for $s\in S$ provides a density, from which the next action/control can be sampled. The usual target function in RL is
\begin{align*}
    \max_{\mu} \E_{\tau}\left[\sum_k \gamma^k r(s_k,a_k)\right], \quad 0<\gamma \leq 1,
\end{align*}
where $\E_{\tau}$ represents the expected value over all possible trajectories $\tau = [s_0^\top,a_0^\top,\dots,s_N^\top]^\top$ of length $N\in \N$. In this manuscript, we make use of an extension of this expression by a trust region idea, where the Kullback-Leibler divergence measures the difference of policies. Thereby, the trust region is not handled as hard constraint, but as penalty term. We refer to~\cite{Schulman.2015,Schulman.2017} for details.

Finally, the actual control $a_k$ given by RL for the current state $s_k\in S$ is 
sampled from the parameterized policy. Hence, the RL feedback law $u_{\text{RL}}$ is given as
\begin{align}\label{eq:uRL}
    u_{\text{RL}}(t) = a_k \sim \pi_{\mu}(\cdot,s_k) \qquad\text{for $t \in [t_k,t_{k+1})$, $k \in \N$.}
\end{align}
We stress that, usually, the policy and the Markov Decision Process are tailored to time-discrete control systems. In this manuscript, the RL control $u_{\text{RL}}$ is assumed to be a sampled-data control 
as it can be seen in (\ref{eq:uRL}). In this way, the discrete-time and continuous-time frameworks of RL and the funnel controller can be combined.

\subsection{Safeguarding feedback law} \label{Sec:FeedbackLaw}
In this section we present the second controller component, namely the so-called \textit{funnel controller}. 
This is a high-gain adaptive controller, which achieves output reference tracking within prescribed error bounds~\eqref{eq:ErrorPerformance}. 
Since this controller is model-free, it safeguards the learning process by compensating undesired control effects from the RL component.
For a given funnel function~$\varphi: \Rp \to \R_{> 0}$ belonging to the set
\begin{equation*}
    \Phi := \setdef{\vp \in W^{1,\infty}(\Rp,\R)}{\inf_{s \ge 0} \vp(s) > 0},
\end{equation*}
and given reference $y_{\rm ref} \in W^{2,\infty}(\Rp,\R^m)$, we formally introduce the funnel control feedback.
First, in virtue of~\cite{berger2021funnel}, we introduce the following auxiliary variables
 \begin{equation} \label{eq:AuxiliaryVariables}
     \begin{aligned}
         e(t) &:= \begin{bmatrix}
         q_1(t) \\ q_2(t)
     \end{bmatrix}- \begin{bmatrix}
         q_{1,\text{ref}}(t) \\ q_{2,\text{ref}}(t)
     \end{bmatrix}, 
     e_1(t) := \varphi(t) e(t), \\
     e_2(t) &:= \varphi(t) \dot{e}(t)+\frac{1}{1-\lVert e_1(t)\rVert^2} e_1(t) \\
     &= \vp(t) \left( \dot e(t)+\frac{1}{1-\lVert e_1(t)\rVert^2} e(t) \right).
     \end{aligned}
 \end{equation}
with $e(t) = e(t,y)$, $e_1(t) = e_1(t,y)$, and $e_2(t) = e_2(t,y,\dot{y})$.
 The funnel control feedback law then reads
 \begin{equation} \label{eq:u_FC}
     u_{\text{funnel}}(t) = -\frac{1}{1-\lVert e_2(t)\rVert^2} e_2(t) .
 \end{equation}
This control law is of adaptive high-gain type.
If the tracking error is small, then the distance to the error boundary~$1/\vp(t)$ is large and no/little input action is required to achieve the tracking task.
If, however, the tracking error approaches its tolerance, then the denominator $1-\vp(t)^2\|e(t)\|^2$ becomes small, and hence, the fraction becomes large. Thus, the error is ``pushed away'' from the funnel boundary and satisfies~\eqref{eq:ErrorPerformance} for all times.
We highlight that, although in~\eqref{eq:u_FC} a possible pole is introduced, it has been proven that the control input is finite for all times, i.e., the error never touches the funnel boundary, cf.
~\cite{berger2021funnel} for arbitrary relative degree.

\subsection{Combined controller}
To achieve the control objective introduced in \Cref{Sec:ControlObjective}, we combine the two controller components.
One first approach could be to set
\begin{align*} 
    u := u_{\text{funnel}}+u_{\text{RL}}.
\end{align*}
However, in this case the funnel controller~\eqref{eq:u_FC} would be active whenever~$e_2(t) \neq 0$ and thus, evaluation of the effectiveness of~$u_{\rm RL}$ is not possible, since $u_{\rm funnel}$ continuously intervenes --~even if RL provides a potentially better control signal due to, e.g., its prediction capabilities.
Therefore, we divide the interior of the funnel into a safe and a safety-critical region, cf. \Cref{Fig:SafeAndSafetycritical}.
More precisely, we introduce an activation threshold~$\lambda \in [0,1)$ to divide the (half-open) interval $[0,1)$ into a safe and a safety-critical region and evaluate the auxiliary variable~$e_2(t)$ w.r.t.\ the activation threshold, cf.~\cite{schmitz2023safe,lanza2023control}.
For a given~$\lambda \in (0,1)$ we define the activation function $\tilde \alpha: [0,1) \to [0,1)$ by
\begin{equation*}
    \tilde \alpha(s) = \max\{0, \|e_2(t)\| - \lambda\}
\end{equation*}
With this, we could define the overall controller
\begin{equation*}
    u := \tilde \alpha(\|e_2(t)\|) \cdot u_{\text{funnel}}+u_{\text{RL}},
\end{equation*}
where the funnel controller is only active, if the auxiliary signal~$e_2$ exceeds the activation threshold~$\lambda$. 
However, this controller is likely to lead to chattering behavior, since whenever~$e_2$ touches the activation threshold, the funnel controller would react and may achieve~$\|e_2\|< \lambda$ immediately.
To avoid possible chattering in the control signal, we introduce a dwell-time activation scheme for the funnel controller as follows 
\begin{equation*}
    \alpha_{t_d}(t,e_2(t)) = \max\left\{ 0, \max_{s \in [t-t_d,t]} \| e_2(s) \| - \lambda \right\}
\end{equation*}
where $t_d > 0$ is a preassigned dwell time.
The previous function ``records'' if the auxiliary variable has exceeded or is above the activation threshold within the past interval of length~$t_d$.
With this activation function, we define the overall control law
\begin{equation} \label{eq:u_total}
    u :=  \alpha_{t_d}(t,e_2(t)) \cdot u_{\text{funnel}}+u_{\text{RL}}.
\end{equation}
Based on the previous considerations, we 
formulate the following feasibility result for the two-component controller.
\begin{thm}
    Consider system~\eqref{eq:Model_input_output_form}. Let a reference trajectory $y_{\rm ref} \in W^{2,\infty}(\Rp,\R^{9})$ and performance bound~$\vp \in \Phi$ be given.
    If the auxiliary variables~\eqref{eq:AuxiliaryVariables} satisfy the initial conditions $\vp(0) \|e_1(0)\| < 1$ and $\vp(0) \|e_2(0)\| < 1$, 
    then the proposed controller~\eqref{eq:u_total} (given \eqref{eq:uRL} and \eqref{eq:u_FC}) achieves the control objective~\eqref{eq:ErrorPerformance}.
\end{thm}
\begin{pf}
    First we observe that the control~$u_{\rm RL}$ from~\eqref{eq:uRL} is piecewise constant and in particular bounded.
    To prove boundedness of the funnel control signal~$u_{\rm funnel}$, and its success in keeping the error variables~$e_1, e_2$ bounded away from 1, usually a contradiction argument is used, cf. the proof of~\cite[Thm.~1.9]{berger2021funnel}.
    Since the respective analysis only considers a small neighborhood on the funnel boundary (cf. the red area in \Cref{Fig:SafeAndSafetycritical}), the incorporation of the activation function~$\alpha_{t_d}(\cdot)$ does not jeopardize applicability of the standard arguments. 
\end{pf}

\section{Numerical example} \label{sec:numerics}
We illustrate the above introduced controller~\eqref{eq:u_total} in a docking maneuver. 
The goal is to reach a predefined position with the arm of the satellite. Meanwhile, the satellite body should move as little as possible to enhance precise gripping. This is a nontrivial task since a change in the arm's position always leads to a force acting on the satellite body. Furthermore, we add collision areas such as solar panels or antennas of the target satellite where the rewards become negative. 

\begin{figure}[ht]
    \begin{subfigure}[b]{0.2\textwidth}
	\begin{center}
		\includegraphics[height=3.1cm]{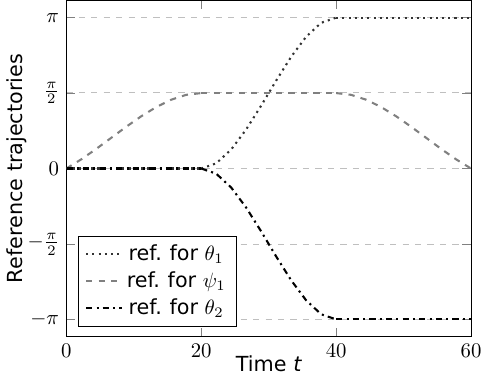}
		\subcaption{Reference trajectory.}
		\label{fig:phi2}
	\end{center}
    \end{subfigure}
    \hspace{0.2cm}
    \begin{subfigure}[b]{0.25\textwidth}
    \begin{center}
		\includegraphics[height=3.5cm]{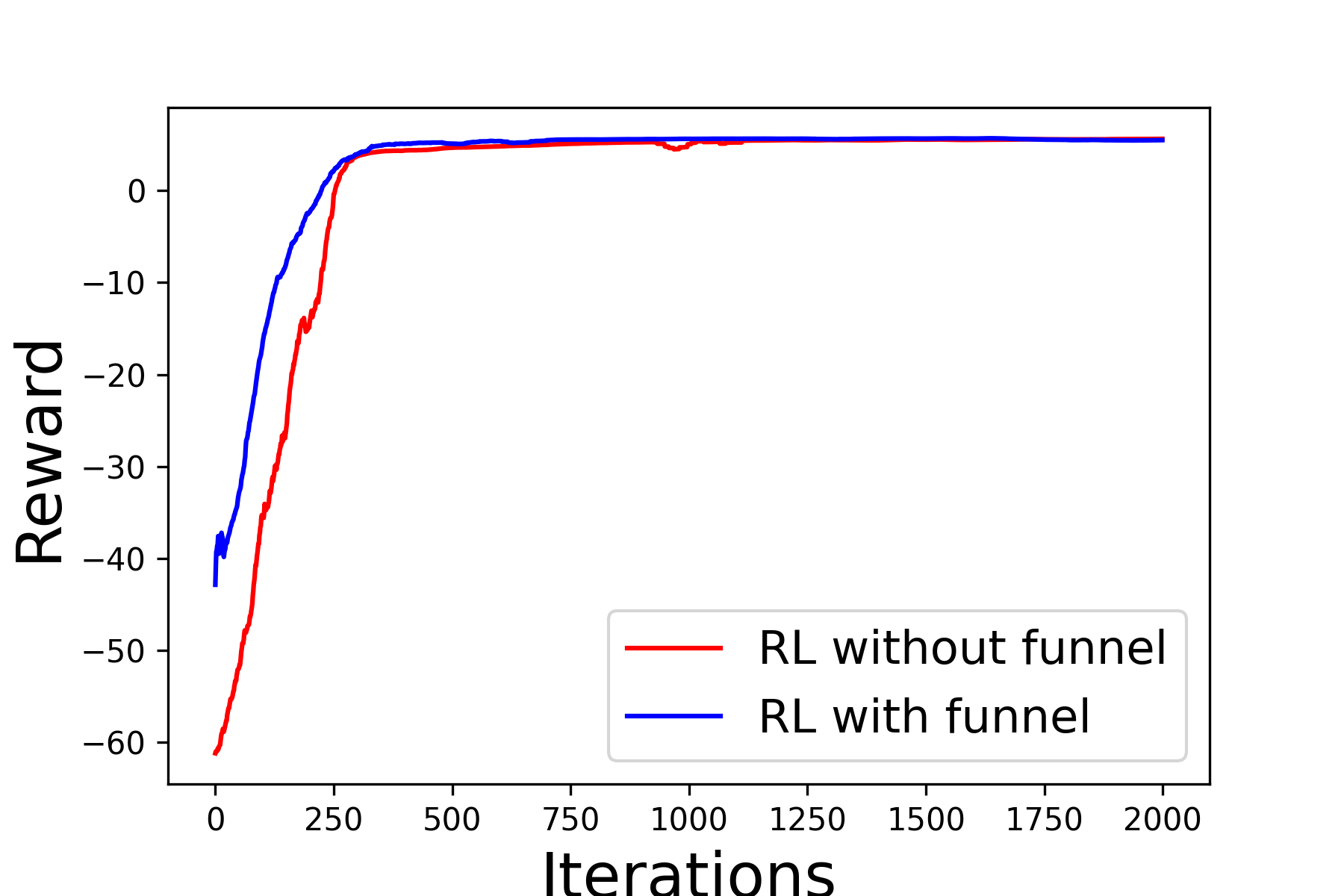}
		\subcaption{Reward during the training.}
		\label{fig:reward_combined}
	\end{center}
    \end{subfigure}
    \caption{Reference trajectory and total reward.}
\end{figure}

With respect to these constraints, we define a reference trajectory, which is safe:
\begin{align*}
    y_{1,\text{ref}}(t) &= 0_{6\times 1}, \forall t \in [0,60], \\
    y_{2,\text{ref}}(t) &= \begin{cases}
        \begin{bmatrix}
            0, & p_1\left(t\right), & 0
        \end{bmatrix}^\top,  \text{if } t\leq 20, \\
        \begin{bmatrix}
            p_2\left(t\right), & \frac{\pi}{2}, & -p_2\left(t\right)
        \end{bmatrix}^\top, \text{if } 20 < t\leq 40,\\
        \begin{bmatrix}
            \pi, & p_3\left(t\right), & -\pi
        \end{bmatrix}^\top,  \text{if } 40 < t,
    \end{cases} \\
    \text{ with } &p_1(t) := \frac{-\pi}{16000}t^3+\frac{\pi}{800}t^2+\frac{\pi}{40}t,\\
    &p_2(t) := \frac{-\pi}{4000} \left(t-20\right)^3+\frac{3\pi}{400} \left(t-20\right)^2,\\
    &p_3(t) := \frac{\pi}{16000} \left(t-40\right)^3-\frac{\pi}{400} \left(t-40\right)^2+\frac{\pi}{2}.
\end{align*}
The reference trajectory $y_{2,\text{ref}}(t)$ is depicted in \Cref{fig:phi2}.
We stress that it is not the goal to just follow the reference trajectory. Instead, the reference trajectory only reveals a collision free trajectory, in whose neighborhood the final trajectory is to be found. For this scenario, we assume the limiting funnel for the deviation from the reference trajectory to be $\vp(t) = 8/\pi$.
It is important to note that within the funnel the positive definiteness of the matrix $g^\top M$ is guaranteed. This means that we can apply the funnel control and it will guarantee that we will not leave the set~$\mathbb{X}$.
Inside this area, we can define the actual optimization task.
For this, the reward function represents the objective function for the PPO algorithm. We stress that this has the potential to cover different control objectives on top of the safeguarding aspect of the funnel controller.

For a better comparability, we run the same RL algorithm without an additional funnel controller. Thus, we distinguish between the reward function for the algorithm, which is supported by a funnel controller, denoted by $r_{\text{funnel}}$, and for a classical pure RL algorithm, denoted by~$r$. For $s_k \in S$ and $a_k \in A$ at time $t_k$, we define:
\begin{multline*}
    r_{\text{funnel}}\left(s_k,a_k\right) \\ = \frac{(1-\lVert e \rVert)}{10}- \int_{t_k}^{t_{k+1}} \alpha_{t_d}(t,e_2) \lVert u_{\text{funnel}}(t) -u_{\rm RL}(t) \rVert \, \text{d}t. 
\end{multline*} 
We point out that the second part of the reward function forces the RL policy to avoid the areas, where funnel control is used, or otherwise to mimic its behavior. In the reward function for the pure RL, the part, which addresses the funnel controller, is omitted. Furthermore, we add a penalty term every time we leave the funnel. The reward function for the pure RL approach reads
\begin{align*}
    r\left(s_k,a_k\right) = \begin{cases}
        \frac{1}{10}(1-\lVert e \rVert), \text{ if } -\frac{1}{\left(\varphi\right)_i} \leq \left(e\right)_i \leq \frac{1}{\left(\varphi\right)_i} \\         \frac{1}{10}(1-\lVert e \rVert)-\frac{60-t}{\Delta h}\left(1+\frac{\lVert e \rVert}{10}\right) , \text{else}.
    \end{cases}
\end{align*}
At this point, we can run the algorithm. For the PPO part of the numerical implementation we use the RL library Rllib (cf. \cite{Liang.2018}). The parameters used for the training can be found in Table~\ref{tb:hyperparametersPPO}. As mentioned before, we apply the algorithm with and without the funnel extension. The reward during the training is depicted in~\Cref{fig:reward_combined}. We observe that in both cases the training process looks promising. The improvements of the rewards are clearly visible. The reward for the approach with funnel control is always greater than the one without funnel control. This is what we expected since high penalty terms for leaving the funnel are avoided by the funnel controller.
\begin{table}
	\begin{center}
		\caption{Algorithm parameters}\label{tb:hyperparametersPPO}
		\begin{tabular}{cc|cc}
			Description & Value & Description &  Value\\\hline
			     RL grid size &$\Delta h=1$ & Discount factor & $\delta=0.9$\\
              Learning rate & lr$=1e^{-4}$ & Activation threshold & $\lambda=0.8$\\
              KL coefficient &  $0.1$ & Dwell time & $t_d= 1$\\
              Batch size &  $128$ &&\\
		\end{tabular}
	\end{center}
\end{table}

\begin{figure}[ht]
    \begin{subfigure}[b]{0.225\textwidth}
	\begin{center}
		\includegraphics[height=3.1cm]{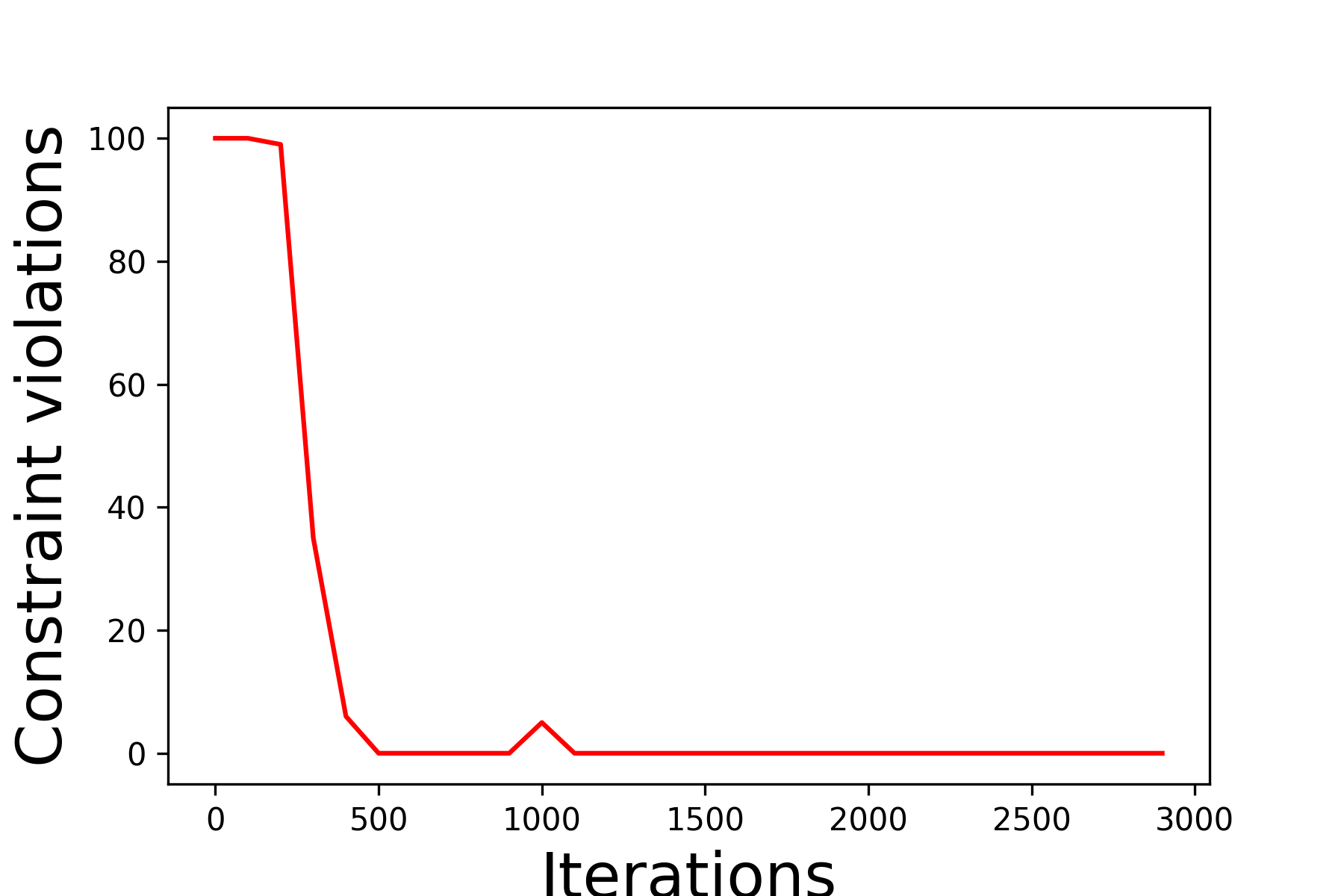}
		\subcaption{Constraint violations without the funnel controller.}
		\label{fig:constraint_violations}
	\end{center}
    \end{subfigure}
    \hspace{0.2cm}
    \begin{subfigure}[b]{0.225\textwidth}
    \begin{center}
		\includegraphics[height=3.1cm]{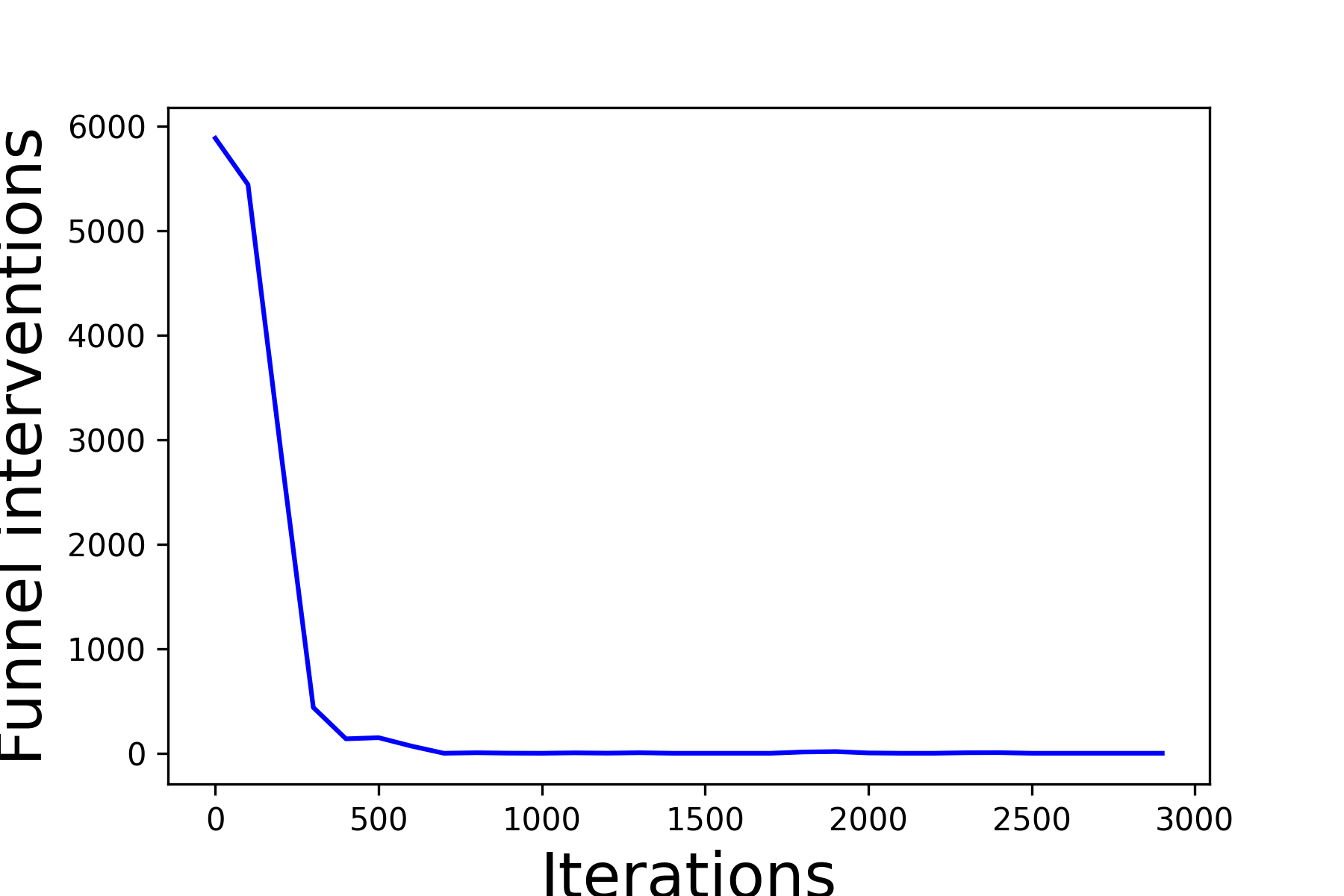}
		\subcaption{Funnel interventions with the funnel controller.}
		\label{fig:funnel_interventions}
	\end{center}
    \end{subfigure}
    \caption{Constraint violations and funnel interventions.}
\end{figure}

In order to show the effect of the funnel controller, we visualized the number of funnel violations in~\Cref{fig:constraint_violations}. It is clear that the use of the funnel controller avoids a violation. 
Without a funnel, the number of violations in $100$ trajectories is high. Especially at the beginning of the training, every trajectory ends with a violation. Keep in mind that a violation could mean a collision of the satellite with its environment. After approximately $350$ iterations, the number of violations decreases and becomes zero at the end of the training since the RL algorithm learns to avoid it.  Nevertheless, we do not have a guarantee that the final policy is safe. 

\begin{figure}[ht]
    \begin{subfigure}[b]{0.225\textwidth}
	\begin{center}
		\includegraphics[height=3.1cm]{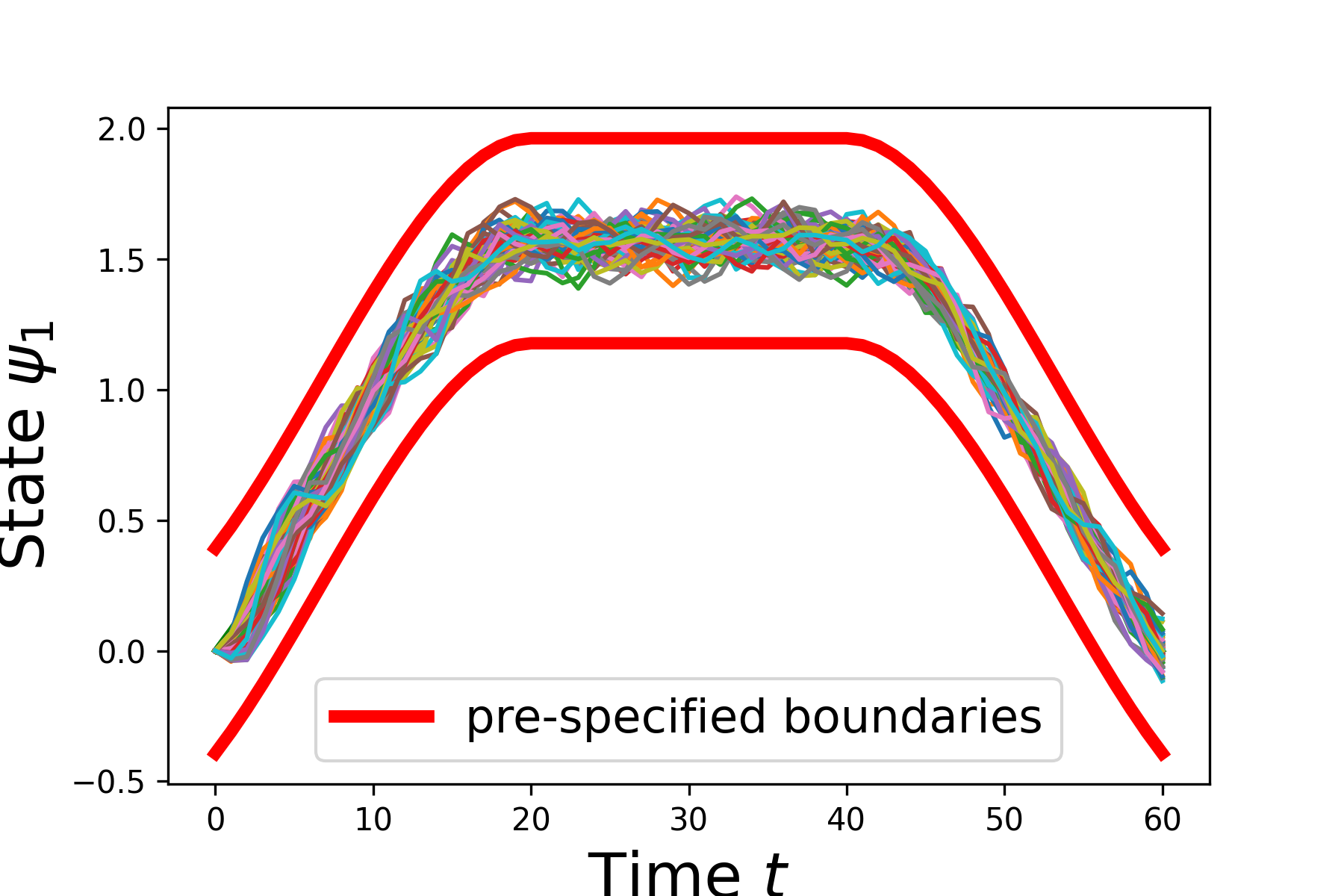}
		\subcaption{Untrained policy.}
		\label{fig:trajectory_before_training}
	\end{center}
    \end{subfigure}
    \hspace{0.2cm}
    \begin{subfigure}[b]{0.225\textwidth}
    \begin{center}
		\includegraphics[height=3.1cm]{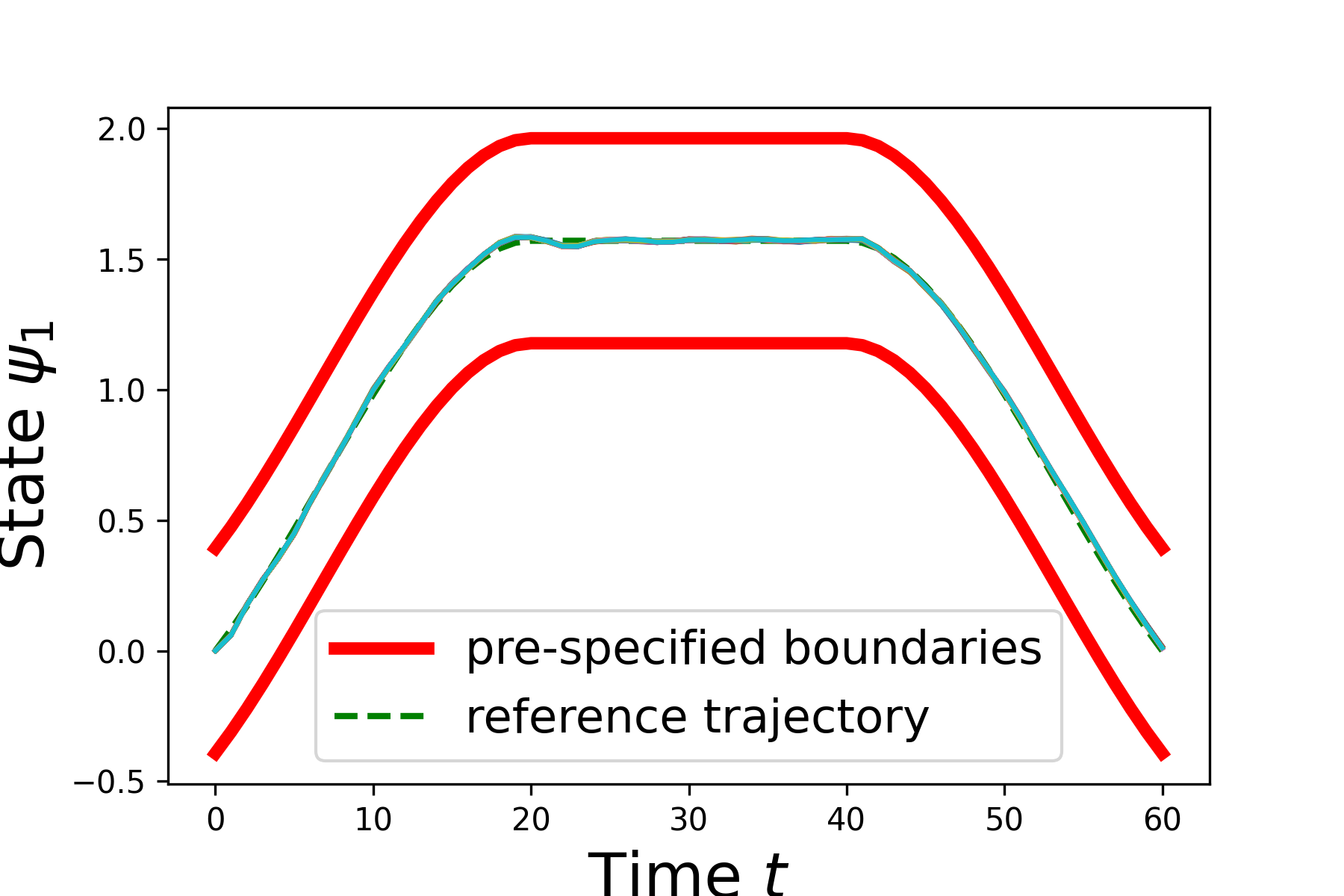}
		\subcaption{Trained policy.}
		\label{fig:trajectory_after_training}
	\end{center}
    \end{subfigure}
    \caption{Generated trajectories for $\psi_1$ with funnel controller.}
\end{figure}

In the funnel-RL case, we are sure by construction that there are no funnel violations. Here, we have a look, how often the funnel controller has to intervene during the training to keep the training safe. In \Cref{fig:funnel_interventions}, we observe that in the beginning, when the policy is not well trained, the funnel control is activated very often. However,
the policy learns how to act such that the funnel controller does not have to intervene. At the end, there was no funnel control intervention needed to generate $100$ trajectories within the desired funnel.
In order to emphasize the effect of the funnel controller, in \Cref{fig:trajectory_before_training}, we plotted the $100$ trajectories of $\psi_1$ for the funnel-RL policy before the training. Thereby, the red thick lines represent the pre-specified boundaries, respectively the funnel. It can be seen that the trajectories are kept inside the funnel. But, in the area close to the reference trajectory the trajectories are chaotic, since the RL part is not trained yet.
For comparison, we plotted $100$ trajectories for the same state $\psi_1$ after the training in~\Cref{fig:trajectory_after_training}. The trajectories are still inside the funnel, but now, all trajectories show a similar behavior due to the trained RL policy.

Overall, we have received a fast controller, which is able to perform the docking maneuver and is protected against collisions with obstacles outside the safety funnel.

\section{Conclusion} \label{sec:conclusion}
To safely conduct docking maneuvers in space, we introduced a two-component controller.
One component is an optimal control policy derived from system data using RL.
To safeguard the learning process and compensate undesired control actions, we add a high-gain
adaptive controller into the control algorithm, which is activated based on evaluations whether the system
is in a safe or a safety-critical region.
Effectiveness of the proposed controller is demonstrated by a numerical example of a satellite docking maneuver with collision avoidance.
In future research, we will focus on the design of an overall control algorithm where the assignment of safe and safety-critical regions accounts for uncertainties both in the model parameters and the trustworthiness of the training state.

\bibliography{ifacconf}             

\end{document}